
 \input amstex
 \documentstyle{amsppt}
 \loadbold
 \magnification=1200
 \baselineskip=15 pt
 \TagsOnRight
 \define\CC{\Cal C}
 \define\CI{\Cal I}
 \define\CM{\Cal M}
 \define\CMs{\Cal M_\sigma}
 \define\CCh{\Cal C_{\text{h}}}
 \define\CCah{\Cal C_{\text{ah}}}
 \define\bCCah{\bar\Cal C_{\text{ah}}}
 \define\CZ{\Cal Z}
 \define\BC{\Bbb C}
 \define\BZ{\Bbb Z}
 \define\bza{\bar z^\ast}
 \define\bX{\bar X}
 \define\bgX{\bar\goth{X}}
 \define\bgZa{\bar\goth{Z}^\ast}
 \define\Uh{\Cal U_h}
 \define\half{{1\over2}}
 \define\thalf{{3\over2}}
 \define\sgn{\text{sgn}}
 \define\eps{\epsilon}
 \define\veps{\varepsilon}
 \define\tr{\text{tr}}
 \define\Mat{\text{Mat}}
 \define\diag{\text{diag}}
 \define\gP{E}
 \define\gX{\goth{X}}
 \define\gY{\goth{Y}}
 \define\gZ{\goth{Z}}
 \define\gV{\goth{V}}
 \define\gN{\goth{N}}
 \define\gM{\goth{M}}
 \define\gQ{\goth{Q}}
 \define\gm{\goth{m}}
 \define\gt{\goth{t}}
 \define\gh{\goth{h}}
 \define\Ad{\Cal A_d} 
 \define\tAd{\tilde\Cal A_d}
 \define\tM{\tilde M}

 \define\bI{\text{\bf I}}
 \define\br{\text{\bf r}}
 \define\bO{\text{\bf 0}}
 %
 \define\ParshalWang{10}
 \define\XxReshetik{9}
 \define\TauberTaft{16}
 \define\Soibelman{13, 14}
 \define\Japanese{1}
 \define\BrezinskiMajid{2}
 \define\Dabrowski{3}
 \define\JSII{8}
 \define\PS{15}
 \define\Drinfeld{4}
 \define\JimboR{6}
 \define\JimboU{5}
 \define\FRT{11}
 \define\JSI{7}
 \define\Rosso{12}
 \topmatter
 \title Antiholomorphic representations for orthogonal and symplectic
 quantum groups$^{\times}$
 \endtitle
 \author Pavel \v S\v tov\'\i\v cek$^{\star\dag}$
 \endauthor
 \leftheadtext{P. \v S\v tov\'\i\v cek}
 \rightheadtext{Quantum orthogonal and symplectic groups}
  \affil
 Arnold Sommerfeld Institute\\
 TU Clausthal, Leibnizstr. 10\\
 D-38678 Clausthal-Zellerfeld, Germany
 \endaffil
 \thanks $^{\star}$ Humboldt Fellow
 \endthanks
 \thanks $^{\dag}$ On leave of absence from:
 Department of Mathematics and Doppler Institute,
 Faculty of Nuclear Science, CTU, Prague, Czech Republic
 \endthanks
 \thanks $^{\times}$ to appear in J. of Algebra
 \endthanks
 \abstract The coadjoint orbits for the series $B_l,\ C_l$ and $D_l$ are
 considered in the case  when the base point is a multiple of a fundamental
 weight. A quantization of the big cell is
 suggested by means of introducing a $\ast$-algebra generated by holomorphic
 coordinate functions. Starting from this algebraic structure the
 irreducible  representations of the deformed universal enveloping algebra
 are derived as acting in the vector space of polynomials in quantum
 coordinate  functions.
 \endabstract
 \endtopmatter
 \document
 \bigpagebreak
 \flushpar {\bf 1. Introduction} 
 \medpagebreak

 Recently a remarkable attention has been paid to quantum homogeneous spaces.
 Particularly interesting is the relation to the representation theory
 in the spirit of the classical method of orbits due to
 Kirillov--Kostant. We quote rather randomly and incompletely
 \cite{\ParshalWang, \XxReshetik, \TauberTaft, \Soibelman, \Japanese,
 \BrezinskiMajid, \Dabrowski} etc.
 Majority of the quoted papers is concerned with a deformation of the
 algebra of holomorphic functions. The quantization of a general coadjoint
 orbit as a complex manifold,
 valid for any compact group from the series $A_l$, $B_l$, $C_l$ and $D_l$,
 has been described in \cite{\JSII} in terms of quantum holomorphic
 coordinate functions living on the big cell. In the same paper, it was shown
 with the help of quantum coherent states that any irreducible
 finite--dimensional representation of the deformed enveloping algebra
 $\Cal U_h$ admits an (anti-) holomorphic realization when acting in a space
 of polynomials in non-commutative variables, the quantum coordinate
 functions.

 The goal of the present paper is to extend the results of the paper
 \cite{\PS}, where the series $A_l$ has been treated, also to orthogonal
 and symplectic groups. For the unitary groups $SU(N)$, the coadjoint orbits
were considered and quantized with the base point being a multiple of a
fundamental weight. These homogeneous spaces coincide exactly with the
Grassmann manifolds. Moreover, antiholomorphic representations of $\Uh$
in the Borel--Weil spirit were derived from the structure of the quantized
orbit. Here we focus on the series $B_l$, $C_l$ and $D_l$, again in the case
when the base point is a multiple of a fundamental weight.
 This means an attempt to suggest a
 quantization of these homogeneous spaces as  real analytic manifolds.
Thus the result should be an algebra $\CC$ generated jointly by the
holomorphic and antiholomorphic quantum coordinate functions $z_{jk}$ and
$z_{st}^{\ \ast}$ living on the big cell. However
the main point consists in deriving representations of $\Uh$ from the algebraic
 structure of $\CC$. Since we restrict ourselves to a special sort of orbits
 the obtained representations are characterized by a lowest weight
 $\lambda$ being an integer multiple of a fundamental weight,
 $-\lambda\in\Bbb Z_+\omega_r$.

 While the general strategy is quite parallel
 to that of the paper \cite{\PS}, the
 commutation relations derived for orthogonal and symplectic groups
 turn out to be considerably more complex and
 less transparent. Apparently this is not the first case when the series
 $B_l$, $C_l$ and $D_l$ are much more difficult to approach than the series
 $A_l$.
 Perhaps one fact standing behind this rather unpleasant observation
 is a property of the R-matrices. While in the $A_l$ case the R-matrix obeys
 the Hecke condition
 $$
 (R_{12}P)^2-(q-q^{-1})R_{12}P-\bI=0,
 $$
 or equivalently,
 $$
 R_{12}-R_{21}^{-1}=(q-q^{-1})P,
 $$
 the R-matrices for the series $B_l$, $C_l$, $D_l$ obey a more complicated
 relation (2.8) below.

Unfortunately, this is why the construction of the algebra $\CC$ was not
clarified entirely and has to some extent speculative character. However
the prescription for the representation of $\Uh$ acting in the vector space
of polynomials in quantum coordinate functions $z_{jk}^\ast$ is derived
quite unambiguously and verified rigorously. In principle, Section 4
concerned with the quantum parameterization of the big cell can be skipped.
But starting directly from the defining relations for the representation
could seem then rather obscure.

 The paper is organized as follows. In Section 2, the notation is introduced
 and a summary of basic relations is given.
Section 3 starts from a quantum version of the orthogonal transformation
of matrices. Furthermore, quantum (anti)-holomorphic coordinate functions
on the big cell are introduced.
A parameterization of the big cell as a real analytic manifold
is suggested in Section 4. The quantization means a specification of
 commutation relations between the coordinate functions.
Section 5 is devoted to verification of the rules according to which
the algebra generated by the quantum antiholomorphic coordinate
functions becomes a left $\Uh$-module. These rules are suggested by the
(partially informal) considerations of Section 4.
The module depends on one free parameter and its values are specified
in Section 6 so that one obtains finite-dimensional
irreducible representations
of $\Uh$. Section 7 contains few final remarks about some points suggested
in Section 4 and not clarified in full detail.

 \bigpagebreak
 \flushpar {\bf 2. Preliminaries}
 \medpagebreak

 The deformation parameter is $q=\text{e}^{-h}$, $h>0$, and the quantum
 numbers are defined as $[x]:=(q^x-q^{-x})/(q-q^{-1})$. The rank of the
 group is denoted by $l$. Let $N:=2l+1$ for the series $B_l$ and
 $N:=2l$ in the case of $D_l$ and $C_l$. The standard basis in $\BC^N$ is
 denoted
 by $\{e_{-l},\dots,e_{-1},e_0,e_1,\dots,e_l\}$ for the series $B_l$ and by
 $\{e_{-l+{1\over2}},\dots,e_{-\half},e_\half,\dots,e_{l-\half}\}$ for the
 series $C_l$ and $D_l$. Set also
 $$
 \align
 &  \eps:=1(:=-1)\ \text{ for the series } B_l
 \text{ and } D_l\ (\text{resp. }C_l), \tag 2.1\\
 & \veps_j:=1\big(:=\sgn(j)\big)
 \text{ for the series } B_l \text{ and } D_l\ (\text{resp. }C_l). \tag 2.2\\
 \endalign$$
 The R-matrix obeying the Yang-Baxter (YB) equation
 $R_{12}R_{13}R_{23}=R_{23}R_{13}R_{12}$ in $\BC^N\otimes\BC^N$ is given by
 \cite{\Drinfeld, \JimboR, \FRT}
 $$
 \align
 & R_{jk,st} :=\delta_{js}\delta_{kt}+(q-q^{\sgn(k-t)})\delta_{jt}\delta_{ks}
 -\veps_k\veps_t
 (q^{-\sgn(k-t)}-q^{-1})q^{-\rho_k+\rho_t+\delta_{k0}\delta_{t0}}
 \delta_{j+k,0}\delta_{s+t,0} \tag 2.3 \\
 & \rho_k :=-k+\eps\,\half\,\sgn\,k.\tag 2.4 \\
 \endalign$$
 Remark: of course, the term $\delta_{k0}\delta_{t0}$ appearing in the
 exponent in the third summand plays a role only for the series $B_l$.

 Further we introduce the $N\times N$ matrices
 $$
 \align
 & \rho:=\text{diag}(\rho_{-(N-1)/2},\rho_{-(N-3)/2},\dots,\rho_{(N-1)/2}),
 \tag 2.5 \\
 & C:=C^0q^\rho=q^{-\rho}C^0,\
 \text{where }\ C^0_{jk}:=\veps_j\delta_{j+k,0}. \tag 2.6 \\
 \endalign$$
 Thus it holds true that $C^t=C^\ast=\eps C^0q^{-\rho}=\eps q^\rho C^0$,
 $C^{-1}=\eps C$ and $C^tC=q^{2\rho}$. Let $K$ be a $N^2\times N^2$ matrix
 with the entries
 $$
 K_{jk,st}:=C^t_{jk}C^{-1}_{\,st}.
 \tag 2.7 $$
 The symbol $P\equiv P_{12}=P_{21}$ stands for the flip operator.

 A survey follows giving the basic relations valid for the matrices
 $R,\ K,\ P$ and C. Quite essential is the first one,
 $$
 R_{12}-R_{21}^{-1}=(q-q^{-1})(P-K_{12}).
 \tag 2.8 $$
 Furthermore,
 $$
 \align
 & R_{21}\equiv PR_{12}P=R_{12}^t,\ R_q^{-1}=R_{q^{-1}}, \tag 2.9 \\
 & C_2^{-1}R_{12}C_2=(R_{12}^{-1})^{t_2},\
 (C_1^{-1})^tR_{12}C_1^t=(R_{12}^{-1})^{t_1}, \tag 2.10 \\
 & C_1C_2R_{12}C_1C_2=R_{21}, \tag 2.11 \\
 & C_1C_2K_{12}=C_1^tC_2^tK_{12}=PK_{12}, \tag 2.12 \\
 & K_{21}\equiv PK_{12}P=K_{12}^t,\ K_{q^{-1}}=K_q^t. \tag 2.13 \\
 \endalign$$
 One can also show that
 $$
 RPK=\eps\,q^{-N+\eps}K
 \tag 2.14 $$
 and hence $R_{12}^{-1}K_{12}=K_{21}R_{21}^{-1}=
 \eps\,q^{N-\eps}PK_{12}$.
 For any $N^2\times N^2$ respectively $N\times N$ matrix $X$ (the dimension
 is clear from the context) we have
 $$
 \align
 & K_{12}X_{12}K_{12}=\tr(K_{12}X_{12})K_{12}, \tag 2.15 \\
 & \tr(K_{12}PX_1)=\eps\,\tr(q^{2\rho}X). \tag 2.16 \\
 \endalign $$
 In the last relation, the trace on the LHS is taken in $\BC^N\otimes\BC^N$
 while on the RHS it is taken in $\BC^N$. Consequently,
 $$
 K_{12}PX_1K_{12}=\eps\,\tr(q^{2\rho}X)K_{12}. \tag 2.17
 $$

 It holds  true that
 $$
 K_{12}X_1=K_{12}C_2^{-1}X_2^tC_2,\
 K_{12}X_2=K_{12}(C_1^{-1})^tX_1^tC_1^t, \tag 2.18
 $$
 particularly,
 $$
 K_{12}R_{31}^{-1}=K_{12}R_{32},\
 K_{12}R_{23}^{-1}=K_{12}R_{13}. \tag 2.19
 $$
 Notice also that
 $$
 K_{12}K_{32}=K_{12}P_{13},\ K_{12}K_{13}=K_{12}P_{23}.
 $$

 Next let us introduce a family of orthogonal projectors in $\BC^N$. The
 symbol $\gP(s)$ stands for the matrix (diagonal) expressed
 in the basis $\{e_k\}$ and corresponding
 to the projector onto $\text{span}\{e_j;\ j\le s\}$; particularly
 $\gP(-(N+1)/2)=0$, $\gP((N-1)/2)=\bI$. We have
 $$
 \align
 & K_{12}\gP(s)_1=K_{12}(\bI-\gP(-s-1)_2), \tag 2.20 \\
 & s+t<0\Longrightarrow K_{12}\gP(s)_1\gP(t)_2=0. \tag 2.21 \\
 \endalign $$
 Furthermore, for all $s$,
 $$
 \gP(s)_1R_{12}(\bI-\gP(s)_1)=0,\
 (\bI-\gP(s)_2)R_{12}\gP(s)_2=0. \tag 2.22
 $$
 One can derive other relations like
 $R_{21}(\bI-\gP(s)_2)=(\bI-\gP(s)_2)R_{21}(\bI-\gP(s)_2)$, etc. Let us fix for
the
 rest of the paper an index $r>0$ ($r\in\{1,\dots,l\}$, in the $B_l$ case,
 and $r\in\{\half,\dots,l-\half\}$, for the series $C_l$ and $D_l$).
 One can show that
 $$
 (\bI-\gP(-r)_1)R_{12}\gP(-r)_1\gP(-r)_2=0. \tag 2.23
 $$
 It is straightforward to verify that
 $$
 \align
 & R_{12}^{-1}\gP(-r)_2K_{12}=\eps\,q^{-2r+1+\eps}\gP(-r)_2PK_{12}\, ,
 \tag 2.24 \\
 & R_{21}(\bI-\gP(-r)_2) K_{12}=\eps\,q^{-2r+1+\eps}(\bI-\gP(-r)_2)PK_{12}
 -(q-q^{-1})\gP(-r)_1PK_{12}. \tag 2.25 \\
 \endalign $$
 The last relation implies
 $$
 \align
 & (\bI-\gP(-r)_1)(\bI-\gP(-r)_2)R_{21}(\bI-\gP(-r)_1)(\bI-\gP(-r)_2)K_{12} \\
 & =\eps\,q^{-2r+1+\eps}(\bI-\gP(-r)_1)(\bI-\gP(-r)_2)PK_{12}. \tag 2.26 \\
 \endalign $$

 The goal of this paper is to determine a relation between the coadjoint
orbits and the representations of the deformed enveloping algebra $\Uh$.
 Let $H_j,\ X_j^\pm$, $j=1,\dots,l$, be the Chevalley generators of $\Uh$
 for a simple Lie algebra of rank $l$. The enumeration of roots in the
 Dynkin diagram is chosen so that the first root is the shortest one or the
 longest one, in the case of $B_l$ or $C_l$, respectively, and the nodes
 1 and 2 are connected with the node 3 in the $D_l$ case.

 The commutation relations defining $\Uh$ can be found in various papers
 \cite{\JimboU, \Drinfeld, \FRT}. $\Uh$ is isomorphic as a $\ast$-algebra to
 the algebra $\Ad$ of quantum functions living on the dual solvable group
 \cite{\FRT}. The generators of $\Ad$ can be arranged into an upper-triangular
 $N\times N$ matrix $\Lambda=(\alpha_{jk})$ with positive (in a convenient
 sense) entries on the diagonal, $\alpha_{jj}>0$, and obeying
 $\prod\alpha_{jj}=1$. The commutation relations can be written as follows
 \cite{\FRT,\JSI}
 $$
 R\Lambda_1\Lambda_2=\Lambda_2\Lambda_1R,\
 \Lambda_1^\ast R^{-1}\Lambda_2=\Lambda_2R^{-1}\Lambda_1^\ast, \tag 2.27
 $$
 jointly with the orthogonality condition,
 $$
 C\Lambda^tC^{-1}=C^{-1}\Lambda^tC=\Lambda^{-1}. \tag 2.28
 $$
 The isomorphism between $\Uh$ and $\Ad$ is given explicitly in terms of
 generators \cite{\FRT},
 $$\align
 & q^{H_j}=(\alpha_{j-1,j-1})^{-1}\alpha_{jj},\quad
 j=1,\dots,l, \\
 & (q-q^{-1})X_j^-=q^{-1/2}(\alpha_{j-1,j-1}\alpha_{jj})^{-1/2}
 \alpha_{j-1,j},\quad j=1,\dots,l, \tag 2.29a \\
 \endalign $$
 for the $B_l$ case ($\alpha_{00}=1$);
 $$\align
 & q^{H_1}=(\alpha_{\half,\half})^2, \\
 & q^{H_j}=(\alpha_{j-\thalf,j-\thalf})^{-1}\alpha_{j-\half,j-\half},
 \quad j=2,\dots,l, \\
 & (q-q^{-1})X_1^-=(q+q^{-1})^{-1/2}q^{-1}\alpha_{-\half,\half}, \\
 & (q-q^{-1})X_j^-=q^{-1/2}
 (\alpha_{j-\thalf,j-\thalf}\alpha_{j-\half,j-\half})^{-1/2}
 \alpha_{j-\thalf,j-\half},\quad j=2,\dots,l, \tag 2.29b \\
 \endalign $$
 for the $C_l$ case;
 $$\align
 & q^{H_1}=\alpha_{\half,\half}\alpha_{\thalf,\thalf}, \\
 & q^{H_j}=(\alpha_{j-\thalf,j-\thalf})^{-1}\alpha_{j-\half,j-\half},
 \quad j=2,\dots,l, \\
 & (q-q^{-1})X_1^-=q^{-1/2}\alpha_{\half,\half}^{\quad 1/2}
 \alpha_{\thalf,\thalf}^{\quad-1/2}\alpha_{-\half,\thalf}, \\
 & (q-q^{-1})X_j^-=q^{-1/2}
 (\alpha_{j-\thalf,j-\thalf}\alpha_{j-\half,j-\half})^{-1/2}
 \alpha_{j-\thalf,j-\half},\quad j=2,\dots,l, \tag 2.29c \\
 \endalign $$
 for the $D_l$ case (let us note that $\alpha_{-\half,\half}=0$); and
 $$
 X_j^+=(X_j^-)^\ast,\quad\forall j, \tag 2.30
 $$
 in all three cases.

 Instead of $\Ad$ we shall deal with the $\ast$-algebra $\tAd$ whose
 generators are arranged into a "positive" matrix
 $$
 M:=\Lambda^\ast\Lambda;\ \text{ hence } M^\ast=M. \tag 2.31
 $$
 The generators of the algebra $\Ad$ can be recovered by decomposing $M$ into
 a product of a lower triangular matrix times an upper triangular matrix.
 With more details, this relationship has bee described in \cite{\PS}.
 From the commutation relations (2.27) one can deduce that
 $$
 M_2R_{12}^{-1}M_1R_{21}^{-1}=R_{12}^{-1}M_1R_{21}^{-1}M_2. \tag 2.32
 $$
 Concerning the orthogonality condition, notice first that (2.27), (2.10) and
 (2.28) imply $C_2\Lambda_1^\ast\Lambda_2^t C_2^{-1}=
 \Lambda_1^\ast\Lambda_2^{-1}=R\Lambda_2^{-1}\Lambda_1^\ast R^{-1}$.
 Consequently, using (2.14),
 $$
 \align
 M_{jk} & =
 \sum_s(C_2^{-1}R_{12}\Lambda_2^{-1}\Lambda_1^\ast R_{12}^{-1}C_2)_{jk,ss}\\
 & =\eps\, q^{N-\eps}\sum_{s}
 (C_2^{-1}R_{12}\Lambda_2^{-1}\Lambda_1^\ast C_2)_{jk,ss}\\
 & =\eps\, q^{N-\eps}\sum_{s}
 [C_2^{-1}R_{12}\Lambda_2^{-1}(C_2^{-1}\Lambda_2^\ast C_2)^tC_2^t]_{jk,ss}\\
 & = \eps\,q^{N-\eps}\sum_{s}(C_2^{-1}R_{12}M_2^{-1}C_2^t)_{jk,ss}.\\
 \endalign $$
 This equality can be rewritten as
 $$
 M_1K_{12}=R_{12}M_2^{-1}R_{12}^{-1}K_{12}
 =\eps\,q^{N-\eps}R_{12}M_2^{-1}PK_{12}. \tag 2.33
 $$
 Summing up, we have
 \proclaim{Definition 2.1} The algebra $\tAd$ is generated by the entries of
the $N\times N$ matrix $M$ and is determined by the relations (2.32), (2.33).
The $\ast$-involution in $\tAd$ is given by $M^\ast=M$
($(M^\ast)_{jk}:=(M_{kj})^\ast$).
 \endproclaim

 Let us show that
 \proclaim\nofrills{} the element $\tr(q^{2\rho}M)$ belongs to the centre of
 $\tAd$, i.e.,
 $$
 \tr(q^{2\rho}M)M=M\tr(q^{2\rho}M). \tag 2.34
 $$
 \endproclaim

 Actually, (2.32) implies
 $(R_{12}^{-1})^{t_2}M_2^tM_1=(R_{12}^{-1}M_1R_{21}^{-1}M_2R_{21})^{t_2}$.
 Owing to (2.10) and (2.19), one gets
 $$\align
 K_{12}M_3^tM_1 & = K_{12}C_3^{-1}R_{23}^{-1}C_3
 (R_{23}^{-1}M_2R_{32}^{-1}M_3R_{32})^{t_3} \\
 & = K_{12}
 (R_{23}^{-1}M_2R_{32}^{-1}M_3R_{32}C_3^tR_{13}^{t_3}(C_3^t)^{-1})^{\ t_3}, \\
 \endalign$$
 and hence
 $$
 K_{12}M_3M_2=K_{12}R_{23}^{-1}M_2R_{32}^{-1}M_3R_{32}R_{13}^{-1}.
 $$
 Having multiplied this equation by $K_{34}P_{34}$ from the left and by
 $K_{34}$ from the right, one can use (2.17) and repeatedly (2.19) to get
 $$\align
 \eps\,\tr(q^{2\rho}M)K_{12}M_2K_{34} & =
 K_{12}K_{34}P_{34}R_{23}^{-1}M_2R_{32}^{-1}M_3R_{32}R_{13}^{-1}K_{34} \\
 & = K_{12}K_{34}P_{34}R_{23}^{-1}M_2R_{32}^{-1}M_3R_{14}R_{42}^{-1}K_{34} \\
 & = K_{12}K_{34}R_{23}^{-1}R_{24}^{-1}P_{34}M_2R_{32}^{-1}R_{42}^{-1}M_3
 K_{34}. \\
 \endalign$$
 Since $K_{34}R_{23}^{-1}R_{24}^{-1}=K_{34}$,
 $K_{34}R_{42}^{-1}R_{32}^{-1}=K_{34}$, we obtain
 $$
 \tr(q^{2\rho}M)K_{12}M_2K_{34}=K_{12}M_2\tr(q^{2\rho}M)K_{34}.
 $$

 Actually much more is known. The elements $\tr(q^{2\rho}M^j)$,
 $j=1,\dots,l$, generate the center of the algebra $\tAd$ \cite{\FRT}.

 An immediate consequence of (2.33) is that
 $$
 \tr(q^{2\rho}M)=\tr(q^{2\rho}M^{-1}). \tag 2.35
 $$
 It is enough to apply $K_{12}P$ from the left and to take into account
 (2.17) and (2.14). Notice also that it holds
 $$
 \tr(q^{2\rho})=\eps+[N-\eps]. \tag 2.36
 $$
 This equality follows also from (2.17), with $X=\bI$, and from (2.8) and
 (2.14) as
 $$
 (q-q^{-1})K_{12}PK_{12}=-K_{12}P(R_{12}-R_{21}^{-1})+(q-q^{-1})K_{12}.
 $$

 \bigpagebreak
 \flushpar {\bf 3. Quantum antiholomorphic coordinate functions
 on the big cell}
 \medpagebreak

 Let us start from an observation which is easy to verify and
 is quite essential for the rest of the paper.
 \proclaim{Lemma 3.1} To any $N\times N$ matrix $X$ there exists a matrix $Y$
 of the same dimension such that
 $$
 K_{12}X_2PR_{12}=K_{12}Y_2. \tag 3.1
 $$
 The same matrix $Y$ obeys also
 $$
 Y_2K_{21}=R_{21}PX_2K_{21}.
 $$
 Moreover, $X^\ast=X$ implies $Y^\ast=Y$.
 \endproclaim
 Thus the entries of $Y$ depend linearly on $X$, namely
 $$
 Y_{jk}=\eps\,\veps_j\veps_k\,q^{-\rho_j-\rho_k-\delta_{jk}}X_{-k,-j}+
 (q-q^{-\sgn(j+k)})q^{-\delta_{j0}\delta_{k0}}X_{jk}-
 \eps\,(q-q^{-1})\delta_{jk}\sum_{s>-j}q^{2\rho_s}X_{ss}. \tag 3.2
 $$
 The term $\delta_{j0}\delta_{k0}$ appearing in an exponent can be again
 omitted for the series $C_l$ and $D_l$. The relation (3.1) can be inverted,
 $$
 K_{12}X_2=K_{12}Y_2PR_{21}^{-1},
 $$
 and it holds true that $X_2K_{21}=R_{12}^{-1}PY_2K_{21}$.

 As a first application, let us show that the orthogonality condition (2.33)
 follows partially from (2.32). By exchanging the indices, one can rewrite
 (2.32) as $M_2R_{12}^{-1}M_1R_{12}=R_{21}M_1R_{21}^{-1}M_2$. Subtracting
 these two equations and using (2.8) we get
 $$
 M_2R_{12}^{-1}M_1K_{12}=K_{21}M_1R_{21}^{-1}M_2. \tag 3.3
 $$
 Let $\tM$ be the matrix defined by
 $$
 K_{12}\tM_2:=K_{12}M_2PR_{21}^{-1}=K_{12}M_2R_{12}^{-1}P. \tag 3.4
 $$
 Multiply (3.4) by $K_{12}$ from the left. Notice that, by (2.17),
 $K_{12}K_{21}=K_{12}PK_{12}P=\eps\,\tr(q^{2\rho})K_{12}P$. Taking into
 account the definition of $\tM$ and (2.17) also on the LHS, one obtains
 $$\align
 & \tr(q^{2\rho}\tM M)K_{12}=\tr(q^{2\rho})K_{12}\tM_2M_2,\text{ i.e.},
 \tag 3.5 \\
 & \tM M=c\bI,\text{ where }c=\tr(q^{2\rho}\tM M)/\tr(q^{2\rho}). \tag 3.6 \\
 \endalign$$
 Since it is true that $R_{12}^{-1}M_1K_{12}=P\tM_1K_{12}$, the LHS of (3.5)
 also equals $\tr(q^{2\rho}M\tM)K_{12}$. For $\tM^\ast=\tM$,
 it holds that $c^\ast=c$ and
 thus $\tM M=(\tM M)^\ast=M\tM$. Consequently, $cM=M\tM M=Mc$ and $c$ belongs
 to the center of $\tAd$. Whence
 \proclaim\nofrills{} the orthogonality condition (2.33) is equivalent to the
equality $c=\eps\,q^{N-\eps}$.
 \endproclaim

 Recall that an index $r$ has been fixed as a positive integer
 (half-integer) not exceeding $l$ (resp. $l-\half$) for the series $B_l$
 ($C_l$ and $D_l$).
 The classical coadjoint orbit can be identified
 with an orbit of the corresponding simple compact group when
 acting on the set of positive $N\times N$ matrices as
 $(U,M)\mapsto U^\ast MU$. The positive matrix $M$ is also required to obey
 $C^0M^tC^0=M^{-1}$. The base point of the orbit is
 $M_0:=\diag(\xi^{-1},\dots,\xi^{-1},1,\dots,1,\xi,\dots,\xi)$, for some
 $\xi>0$. The eigen-values $\xi$ and $\xi^{-1}$ have the same multiplicity
 $(N-2r+1)/2$. Here we are doing one exception: in the $D_l$ case,
 with $r=3/2$, the orbit corresponds to a multiple of the weight
 $\omega_1+\omega_2$ rather than to a multiple of the fundamental weight
 $\omega_2$ (cf. Sec. 6).

 The algebra $\CCh$ generated by the quantum "holomorphic coordinate
 functions" on the big cell has been introduced in \cite{\JSII}. Let us recall
 the definition. Let $\CZ$ be an $N\times N$ matrix split into blocks as
 $$
 \CZ=
 \pmatrix \bI_1 & Z_1 & Z_2 \\
 \bO & \bI_2 & Z_3 \\
 \bO & \bO & \bI_3  \endpmatrix, \tag 3.7
 $$
 with the diagonal unit blocks $\bI_1$ and $\bI_3$ having the same dimension
 $\half(N-2r+1)\times\half(N-2r+1)$ and with $\bI_2$ having the dimension
 $(2r-1)\times(2r-1)$. Introduce an $N^2\times N^2$ matrix $Q$ by
 $$\align
 Q_{jk,st} & = R_{jk,st},\quad\text{if }\ j,s\le-r\ \text{ or }\ -r<j,s<r\
 \text{ or }\ r\le j,s, \\
 & = 0,\quad\text{otherwise}. \tag 3.8
 \endalign$$
 Notice that after replacing $j,s$ by $k,t$ everywhere in the condition
 following "if" on the first line, we get the same matrix $Q$. The commutation
 relation defining the algebra $\CCh$ reads as
 $$
 R_{12}Q_{12}^{-1}\CZ_1Q_{12}\CZ_2=Q_{21}^{-1}\CZ_2Q_{21}\CZ_1R_{12}. \tag 3.9
 $$
 In addition, $\CZ$ is required to obey the "orthogonality" condition
 $$
 \delta_{jk}=\sum_s(\CZ_2C_2^{-1}Q_{12}\CZ_2^tQ_{12}^{-1}C_2)_{kj,ss}.
 $$
 Since $C_1^{-1}C_2^{-1}Q_{12}=Q_{21}C_1^{-1}C_2^{-1}$ and owing to (2.12),
this
 equality can be rewritten as
 $$
 K_{12}=\CZ_1Q_{12}C_1^{-1}\CZ_1^tC_1Q_{12}^{-1}K_{12}. \tag 3.10
 $$
 The relations adjoint to (3.9) and (3.10) define another algebra $\CCah$
 generated by the quantum "antiholomorphic coordinate functions" $z_{jk}^\ast$.

 The relations (3.9), (3.10) can be simplified. To this end we set
 $$
 Z:=(Z_1,Z_2),\
 \gZ:=\pmatrix \bO & Z \\ \bO & \bO \endpmatrix. \tag 3.11
 $$
 Thus $Z=(z_{jk})$ is a $\half(N-2r+1)\times\half(N+2r-1)$ block with the
 indices restricted by $j\le-r,\ -r<k$, and $\gZ$ is an $N\times N$ matrix; the
 dimensions of the blocks are determined implicitly. Set
 $$
 \gP^-:=\gP(-r),\ \gP^0:=\gP(r-1)-\gP(-r),\ \gP^+:=\bI-\gP(r-1). \tag 3.12
 $$
 Thus $\bI=\gP^-+\gP^0+\gP^+$ is a resolution of unity. We also have
 $\gZ=\gP^-\gZ(\bI-\gP^-)$. The matrix $Q$ can be written as (c.f. (2.22))
 $$\align
 Q_{12} & = \gP^-_2R_{12}\gP^-_2+\gP^0_2R_{12}\gP^0_2+\gP^+_2R_{12}\gP^+_2 \\
 & = R_{12}\gP^-_2+\gP^0_2R_{12}\gP^0_2+\gP^+_2R_{12}. \tag 3.13 \\
 \endalign$$
 A similar equality holds true also provided the projectors
 $\gP^-$, $\gP^0$, $\gP^+$ are
 applied in the first factor of the tensor product rather than in the second
 one. Particularly,
 $$
 \gP_1^-R_{12}Q_{12}^{-1}=\gP_1^-. \tag 3.14
 $$
 Furthermore, we have $Q_{q^{-1}}=Q_q^{-1}$. Now, let us multiply (3.9) by
 $\gP_1^-\gP_2^-$ from the left. Using (3.13), (3.14), (2.22) and noticing
 that $\gP^-\CZ=\gP^-+\gZ$, one can derive easily
 $$
 (\gP_1^-+\gZ_1)R_{12}(\gP_2^-+\gZ_2)R_{12}^{-1}=
 R_{21}^{-1}(\gP_2^-+\gZ_2)R_{21}(\gP_1^-+\gZ_1). \tag 3.15
 $$
 Similarly, by multiplying (3.10) by $\gP_1^-\gP_2^-$ from the left and using
 (2.21), the equality $\gP_2^-Q_{12}^{\pm1}=R_{12}^{\pm1}\gP_2^-$, (2.24) and
 (2.18), one obtains
 $$
 (\gP_2^-+\gZ_2)R_{21}(\gP_1^-+\gZ_1)K_{12}=0. \tag 3.16
 $$

 In fact, one can show that
 \proclaim{Lemma 3.2} The relation (3.16) follows from (3.15).
 \endproclaim

 Actually, interchange the sides and the indices 1, 2 in (3.15) to get
 $$
 (\gP_1^-+\gZ_1)R_{12}(\gP_2^-+\gZ_2)R_{21}=
 R_{12}(\gP_2^-+\gZ_2)R_{21}(\gP_1^-+\gZ_1).
 $$
 Subtracting this equation from (3.15) and using (2.8), one arrives at
 $$
 K_{21}(\gP_1^-+\gZ_1)R_{12}(\gP_2^-+\gZ_2)=
 (\gP_2^-+\gZ_2)R_{21}(\gP_1^-+\gZ_1)K_{12}.
 $$
 Now it is enough to multiply the last equality by $\gP_1^-\gP_2^-$ from the
 left.

Summing up we have
\proclaim{Definition 3.3} The algebra $\CCh$ is generated by the entries of
matrix $Z=(z_{jk})$ introduced in (3.11) and is determined by the relation
(3.15). The algebra $\CCah$ is generated by the entries of the matrix
$Z^\ast=(z_{kj}^\ast)$ and determined by the adjoint relation to (3.15).
\endproclaim

 Let $Y$ be a matrix defined by the relation
 $$
 Y_1K_{12}:=\eps\,q^{-2r+1+\eps}R_{12}P(\gP_1^-+\gZ_1)K_{12}. \tag 3.17
 $$
 One can deduce from (3.2) that
 $$
 Y=\gP^++\gV,\quad\text{where }\
 \gV:=\pmatrix \bO&V \\ \bO&\bO \endpmatrix \tag 3.18
 $$
 and the ranges of indices in the
 $\half(N+2r-1)\times\half(N-2r+1)$ block $V=(v_{jk})$ are given by
 $j<r,\ r\le k$. One can write explicitly
 $$
 v_{jk}=q^{-2r+1+\eps}\veps_j\veps_k\big(q^{-\rho_j-\rho_k}z_{-k,-j}+
 \eps\,\chi_{j\le-r}
 (q-q^{-\sgn(j+k)})z_{jk}\big), \tag 3.19
 $$
 where $\chi_{j\le-r}=1\ (0)$ provided the condition $j\le-r$ is (is not)
 satisfied. Thus the orthogonality condition (3.16) means that
 \proclaim{Lemma 3.4} The matrix $V=(v_{jk})$ defined by (3.19) fulfills
 $$
 \pmatrix \bI&Z \endpmatrix \pmatrix V\\ \bI\endpmatrix =0. \tag 3.20
 $$
 \endproclaim

\noindent{\it Notation}. The following shorthand notation will be used
mainly in the proofs. Set
$$
\gamma:=q-q^{-1},\quad p:=\eps\, q^{-2r+1+\eps}. \tag 3.21
$$
Furthermore,
$$
\gX:=\gP^-+\gZ^\ast,\quad \gY:=\gP^++\gV^\ast. \tag 3.22
$$
So $\gX_1R_{12}\gX_2R_{21}=R_{12}\gX_2R_{21}\gX_1$,
$K_{21}\gY_1=pK_{21}\gX_1P_{12}R_{21}$ and
$\gY_1K_{12}=pR_{12}P_{12}\gX_1K_{12}$. Notice that $\gX^2=\gX$ and
$\gY^2=\gY$ (since $(\gZ^\ast)^2=0=(\gV^\ast)^2$). Moreover,
$\gY\gX=0$ (Lemma 3.4) and $\gX\gY=0$ (evident). In addition to
the equality $K_{21}\gX_1R_{12}\gX_2=0$ (Lemma 3.2) we have also
$$
\gX_1R_{12}\gX_2K_{21}=0, \tag 3.23
$$
as a consequence of $\gP^-_1\gP^-_2K_{12}=0$ and $\gX=\gX\gP^-$.

\proclaim{Lemma 3.5}
The matrix $\gV$ defined in (3.17), (3.18) obeys
$$
(\gP^+_1+\gV^\ast_1)R_{21}^{-1}(\gP^+_2+\gV^\ast_2)R_{12}^{-1}=
R_{21}^{-1}(\gP^+_2+\gV^\ast_2)R_{12}^{-1}(\gP^+_1+\gV^\ast_1)
\tag 3.24 $$
and consequently
$$
(\gP^+_1+\gV^\ast_1)R_{21}^{-1}(\gP^+_2+\gV^\ast_2)K_{21}=0.
\tag 3.25 $$
Moreover, it holds true that
$$
R_{21}(E^-_1+\gZ^\ast_1)R_{12}(\gP^+_2+\gV^\ast_2)=
(\gP^+_2+\gV^\ast_2)R_{21}(E^-_1+\gZ^\ast_1)R_{12}
\tag 3.26 $$
and also (equivalently)
$$
(E^-_1+\gZ^\ast_1)R_{21}^{-1}(\gP^+_2+\gV^\ast_2)R_{12}^{-1}=
R_{21}^{-1}(\gP^+_2+\gV^\ast_2)R_{12}^{-1}(E^-_1+\gZ^\ast_1).
\tag 3.27 $$
\endproclaim

\demo{Remark}
The algebra $\CCah$ is also generated by the entries of the matrix
$V^\ast=(v_{kj}^\ast)$ and is determined by the relation (3.24).
\enddemo

\demo{Proof}
The equality
$$
K_{13}K_{24}\gX_3R_{34}\gX_4R_{43}=K_{13}K_{24}R_{34}\gX_4R_{43}\gX_3
$$
leads to (using also (2.19) )
$$
K_{13}K_{24}\gY_3P_{13}R_{31}^{-1}\gY_4P_{14}R_{42}^{-1}R_{32}^{-1}R_{43}=
K_{13}K_{24}\gY_4P_{14}R_{42}^{-1}\gY_3P_{13}R_{31}^{-1}R_{41}^{-1}R_{12}.
$$
Using twice the YB equation on the both sides one obtains
$$
K_{13}K_{24}\gY_3\gY_4R_{23}R_{21}=K_{13}K_{24}\gY_4\gY_3R_{14}R_{34}.
$$
The relation (2.19) then gives
$$
K_{13}K_{24}R_{43}\gY_3R_{43}^{-1}\gY_4=
K_{13}K_{24}\gY_4R_{34}^{-1}\gY_3R_{34}.
$$
To show (3.26) one can proceed similarly starting from
$$
K_{31}\gX_2R_{21}\gX_1R_{21}^{-1}=K_{31}R_{12}^{-1}\gX_1R_{12}\gX_2.
$$
Inserting $R_{12}=R_{21}^{-1}+\gamma P_{12}-\gamma K_{12}$ resp.
$R_{21}=R_{12}^{-1}+\gamma P_{12}-\gamma K_{21}$ on the LHS resp. the RHS
of the equality
$$
\gX_1R_{12}\gY_2R_{12}^{-1}=R_{21}^{-1}\gY_2R_{21}\gX_1 ,
$$
one obtains easily (3.27). \qed
\enddemo

 \bigpagebreak
 \flushpar {\bf 4. Parameterization of the quantum orbit}
 \medpagebreak

The character of this section is to some extent speculative. However one
can extract from the procedure presented here a prescription for an
action of $\tAd$ on $\CCah$. Its rigorous verification is the goal of
the next section.

 Introduce $N\times N$ matrices $\gQ^\pm$ by
 $$
 \gQ^-:=\pmatrix \bI \\ Z^\ast \endpmatrix (\bI+ZZ^\ast)^{-1}
 \pmatrix \bI&Z \endpmatrix,\
 \gQ^+:=\pmatrix V \\ \bI \endpmatrix (\bI+V^\ast V)^{-1}
 \pmatrix V^\ast & \bI \endpmatrix . \tag 4.1
 $$
 They possess projector-like properties: $(\gQ^\pm)^2=\gQ^\pm$ and
 $(\gQ^\pm)^\ast=\gQ^\pm$. Moreover, in virtue of (3.20),
 $\gQ^-\gQ^+=\gQ^+\gQ^-=0$. Relate to a triplet of parameters
 $\xi_0,\xi_+,\xi_-$ the matrix
 $$
 \gM:=\xi_0\bI+(\xi_+-\xi_0)\gQ^++(\xi_- -\xi_0)\gQ^-. \tag 4.2
 $$
 \proclaim{Definition 4.1}
 Denote by $\CC$ the $\ast$-algebra generated by $z_{jk},z_{jk}^\ast$ and
 determined by (3.15),
 $$
 \pmatrix \bI&Z \endpmatrix_2R_{12}^{-1}\gM_1R_{21}^{-1}
 \pmatrix -Z\\ \bI\endpmatrix_2 =0, \tag 4.3
 $$
 and by
 $$
 \tilde\gM=\eps\,q^{N-\eps}\gM^{-1},\ \text{ where }\
 K_{12}\tilde\gM_2:=K_{12}\gM_2R_{12}^{-1}P. \tag 4.4
 $$
 Finally, the element
 $$
 \gm:=\tr(q^{2\rho}\gM) \tag 4.5
 $$
 is required to be central in $\tAd$. Of course, the adjoint relations should
 be fulfilled as well.
 \endproclaim

\demo{Remark}
The structure of the relation (4.3) is not transparent enough and so it is
not clear whether $\CC\not=0$, i.e., whether the unit is not contained in
the ideal generated by the defining relations. Nevertheless, in this section
we shall proceed optimistically as if $1\in\CC$.
\enddemo

 The equality (4.3) means exactly that
 $$
 \gQ_2^-R_{12}^{-1}\gM_1R_{21}^{-1}(\bI-\gQ_2^-)=0. \tag 4.6
 $$
 From the Hermitian property of $\gQ^-$ and $\gM$ and from (2.9) one deduces
that
 this is the same as
 $$
 \gQ_2^-R_{12}^{-1}\gM_1R_{21}^{-1}=R_{12}^{-1}\gM_1R_{21}^{-1}\gQ_2^-.
 \tag 4.7 $$

 It holds again true that (in the same way as in (2.35))
 $$
 \tr(q^{2\rho}\gM)=\tr(q^{2\rho}\gM^{-1}). \tag 4.8
 $$
 Set
 $$
 \gt_\pm:=\tr(q^{2\rho}\gQ^\pm). \tag 4.9
 $$
 The element $\gm$ is central in $\tAd$ and, by (4.2),
 $$
 \gm=\xi_0\tr(q^{2\rho})+(\xi_+-\xi_0)\gt_++(\xi_--\xi_0)\gt_- . \tag 4.10
 $$
 Now, inserting the expression (4.2) for $\gM$ and the expression
 $$
 \gM^{-1}=\xi_0^{-1}\bI+(\xi_+^{-1}-\xi_0^{-1})\gQ_+ +
 (\xi_-^{-1}-\xi_0^{-1})\gQ_-  \tag 4.11
 $$
 for $\gM^{-1}$ into the equality (4.8) one gets a linear dependence between
 $\gt_+$ and $\gt_-$. It follows that, in the generic case, both $\gt_+$
 and $\gt_-$ are expressible in terms of $\gm$, and thus $\gt_+,\gt_-$ are
 central elements in $\tAd$ as well.

 Furthermore, the substitution of (4.2) and (4.11) for $\gM$ and $\gM^{-1}$,
 respectively, into (4.4) leads to the equality
 $$\align
 & K_{12}(q^{N-\eps}\xi_0^{-1}-\eps\,\xi_0PR_{21}^{-1})+
 K_{12}\gQ_2^+\big(q^{N-\eps}(\xi_+^{-1}-\xi_0^{-1})-
 \eps\,(\xi_+-\xi_0)PR_{21}^{-1}
 \big) \\
 & +K_{12}\gQ_2^-\big(q^{N-\eps}(\xi_-^{-1}-\xi_0^{-1})-
 \eps\,(\xi_--\xi_0)PR_{21}^{-1}
 \big)=0.
 \tag 4.12 \\
 \endalign$$
 Provided the matrix
 $q^{N-\eps}(\xi_+^{-1}-\xi_0^{-1})\bI-\eps\,(\xi_+-\xi_0)PR_{21}^{-1}$
 is regular one
 can express $\gQ^+$ in terms of $\gQ^-$. The solution can be find with the
 help of the following ansatz:
 $$
 K_{12}\gQ_2^+=\mu K_{12}+\eta K_{12}\gQ_2^-+\zeta K_{12}\gQ_2^-PR_{12}.
 \tag 4.13 $$
 Whence,
 $$
 K_{12}\gQ_2^+PR_{21}^{-1}=\eps\,q^{N-\eps}\mu K_{12}+\eta
K_{12}\gQ_2^-PR_{21}^{-1}
 +\zeta K_{12}\gQ_2^-. \tag 4.14
 $$
 Owing to (2.8),
 $$
 K_{12}\gQ_2^-PR_{12}= K_{12}\gQ_2^-PR_{21}^{-1}+(q-q^{-1})K_{12}\gQ_2^-
 -\eps\,(q-q^{-1})\gt_-K_{12}. \tag 4.15
 $$
 Inserting (4.13), (4.14) into (4.12), taking into account (4.15) and
 comparing the coefficients at the terms $K_{12}$, $K_{12}\gQ_2^-$ and
 $K_{12}\gQ_2^-PR_{21}^{-1}$, one arrives at a system of equations
 $$\align
 & \xi_0^{-1}-\xi_0+[(\xi_+^{-1}-\xi_0^{-1})-(\xi_+-\xi_0)]\mu
 -\eps\,(q-q^{-1})(\xi_+^{-1}-\xi_0^{-1})\gt_-\zeta=0, \tag 4.16a \\
 & \xi_-^{-1}-\xi_0^{-1}+(\xi_+^{-1}-\xi_0^{-1})\eta+
 [-\eps\,q^{N-\eps}(\xi_+-\xi_0)+(q-q^{-1})(\xi_+^{-1}-\xi_-^{-1})]\zeta=0,
 \tag 4.16b\\
 & \xi_--\xi_0+(\xi_+-\xi_0)\eta-\eps\,q^{N-\eps}(\xi_+^{-1}-\xi_0^{-1})\zeta
 =0. \tag 4.16c \\
 \endalign$$
 In the generic case, the unique solution $(\mu,\eta,\zeta)$ determines the
 RHS in (4.13).

 We claim that

 \proclaim\nofrills{}
 there exists a $\ast$-algebra morphism $\psi:\tAd\to\CC$ prescribed by its
 values on the generators, namely $\psi(M)=\gM$.
 \endproclaim

 To verify this statement one has to check the definition of $\tAd$, i.e.,
 (2.32), (2.33), and the definition of $\CC$, i.e., (4.3), (4.4), giving
 also (4.7). The rest is a consequence of (4.13) and the following
 proposition.

 \proclaim{Proposition 4.2} Let $X$ and $\gN$ be $N\times N$ matrices with
 entries from some associative algebra and assume that they obey the equality
 $$
 X_1R_{21}^{-1}\gN_2R_{12}^{-1}=R_{21}^{-1}\gN_2R_{12}^{-1}X_1 .\tag 4.17
 $$
 Then the matrix $Y$ defined by the relation
 $K_{12}Y_2:=K_{12}X_2PR_{12}$ fulfills an analogous equality,
 $$
 Y_1R_{21}^{-1}\gN_2R_{12}^{-1}=R_{21}^{-1}\gN_2R_{12}^{-1}Y_1 .
 $$
 \endproclaim

 \demo{Proof}
 The symbol $\tr_1(\cdot)$ stands for the trace taken in the first factor
 of a tensor product. Notice also that if $A_{12},\ B_{123}$ are two matrices
 obeying $A_{12}B_{345}=B_{345}A_{12}$ (the entries commute) then
 $$\align
 \tr_1(K_{12}A_{23}^{\ t_3}B_{123})^t & =
 \tr_1(B_{123}^{\quad t}A_{23}^{\ t_2}K_{21}) \\
 & = \tr_1(B_{123}^{\quad t}P_{12}C_1C_2A_{23}K_{21}). \tag 4.18 \\
 \endalign $$
 In the last equality we have used (2.18) and (2.12). Furthermore, (4.17)
 can be rewritten as
 $$
 X_1\gN_2^t=(R_{21}^{-1}\gN_2R_{12}^{-1}X_1R_{12})^{t_2}
 [(R_{21}^{-1})^{t_2}]^{-1} .
 $$
 Now, taking into account that $\tr_1\,K_{12}=\bI_2$ and
 applying successively the definition of $Y$; (4.17) and (2.10); (4.18); then
 (2.11), YB equation and (2.19); the definition of $Y$ and YB equation;
 (2.19) and (2.12); and finally (2.18); we get
 $$\align
 (Y_2R_{32}^{-1}\gN_3)^{t_3} & =
 \tr_1(K_{12}X_2P_{12}R_{12}\gN_3^t(R_{32}^{-1})^{t_3}) \\
 & =\tr_1\big(K_{12}(R_{32}^{-1}\gN_3R_{23}^{-1}X_2R_{23})^{t_3}
 C_3^tR_{32}^{-1}P_{12}R_{12}R_{32}(C_3^t)^{-1}\big) \\
 & =\tr_1(C_3^{-1}R_{23}R_{21}R_{13}^{-1}C_1C_2C_3R_{32}^{-1}
 \gN_3R_{23}^{-1}X_2R_{23}K_{21})^t \\
 & =\tr_1(C_1C_2R_{31}^{-1}R_{12}\gN_3R_{23}^{-1}R_{13}^{-1}X_2K_{21})^t \\
 & =\tr_1(C_1C_2R_{31}^{-1}\gN_3R_{13}^{-1}R_{23}^{-1}P_{12}Y_2K_{21})^t \\
 & =\tr_1(C_1R_{31}^{-1}\gN_3R_{13}^{-1}Y_1R_{13}C_1^{-1}K_{21})^t \\
 & =\tr_1\big((R_{32}^{-1}\gN_3R_{23}^{-1}Y_2R_{23})^{t_2}K_{21}\big)^t \\
 & =(R_{32}^{-1}\gN_3R_{23}^{-1}Y_2R_{23})^{t_3}.\quad\qed \\
 \endalign$$
 \enddemo

 The construction of representations is based on the same idea as in the case
 of the series $A_l$ \cite{\PS}. Owing to the morphism
 $\psi:\tAd\to\CC$, $\CC$ is a left $\tAd$-module. Let $\CI$ be the left
 ideal in $\CC$ generated by the "holomorphic" elements $z_{jk}$
 and $\CC/\CI$ be the
 factor module. As a next step, one should identify $\CC/\CI$ with $\CCah$
 as a vector space (recall that $\CCah$ is a unital subalgebra in $\CC$
 generated only by $z_{jk}^\ast$)
 and consider the cyclic submodule $\CM$ with the cyclic
 vector $1\in\CCah$. Though the structure of relations defining the algebra
 $\CC$ has not been clarified in full detail yet  we shall accept the
 identification $\CC/\CI\equiv\CCah$ as a hypothesis.
 It has been also confirmed by explicit calculations in the lowest rank
 cases.
 Nevertheless, further we shall arrive at a quite unambiguous prescription
 for the action of $\tAd$ (and hence of $\Uh$) on $\CM$
 (Proposition 5.4). In what follows, the central
 dot "$\cdot$" stands for this action.

 \proclaim{Lemma 4.3}
 Assuming $\CC/\CI\equiv\CCah$ it holds true that
 $$
 (\bI+ZZ^\ast)^{-1}\cdot1=\eps\,{(1-\mu)\over\zeta} q^{-2r+1+\eps}\bI .
 \tag 4.19 $$
 Consequently,
 $$
 \gt_-\cdot1=\eps\,{(1-\mu)\over\zeta}q^{(N-2r+1)/2}[(N-2r+1)/2] , \tag 4.20
 $$
 ($[x]=(q^x-q^{-x})/(q-q^{-1})$ ).
 \endproclaim

 \demo{Proof}
 Set temporarily $X_{st}:=(\bI+ZZ^\ast)^{-1}_{\ st}\cdot1$. From (4.1),
 (4.13) and (2.3), one finds that for $j,k\ge r$,
 $$\align
 (\bI+V^\ast V)^{-1}_{\ jk}= & \, (\mu \bI+\eta\gQ^-)_{jk}+
 \eps\,\veps_j\veps_k\zeta\gQ^-_{-k,-j}
 q^{-\rho_k-\rho_j}+\zeta(q-q^{-\sgn(j+k)})\gQ^-_{jk} \\
 & -\,\eps\,\zeta\delta_{jk}\sum_\sigma(q^{\sgn(k+\sigma)}-q^{-1})
 q^{2\rho_\sigma} \gQ^-_{\sigma\sigma} . \tag 4.21 \\
 \endalign$$
 Notice that $\gQ^-_{st}\cdot1=0$, for $t>-r$,
 $\gQ^-_{st}=(\bI+ZZ^\ast)^{-1}_{\ st}$, for $s,t\le-r$, and
 $(\bI+V^\ast V)^{-1}\cdot1=\bI$. Thus (4.21) applied on the unit yields
 (we set $s=-k,\ t=-j;\ s,t\le-r$)
 $$
 \delta_{st}=\mu\delta_{st}-\eps\,\veps_s\veps_t\zeta\,q^{\rho_s+\rho_t}X_{st}
 -\eps\,\zeta\,\delta_{st}\sum_{s\le\sigma\le-r}
 (q^{\sgn(\sigma-s)}-q^{-1})q^{2\rho_\sigma}X_{\sigma\sigma}. \tag 4.22
 $$
 It follows that $X_{st}=0$, for $s\not=t$. Set $x_s:=q^{2\rho_s}X_{ss}$.
 One deduces from (4.22) that
 $$
 1=\mu+\eps\,\zeta x_s-\eps\,\zeta\sum_{s\le\sigma\le -r}
 (q^{\sgn(\sigma-s)}-q^{-1})x_s .
 $$
 The last equality amounts to a recurrent relation for $x_s$ leading
 immediately to the solution. \qed
 \enddemo

 As a consequence one finds that
 $$
 \gQ^-\cdot1=\eps\,{1-\mu\over\zeta}\,q^{-2r+1+\eps}
 \pmatrix \bI&\bO\\ Z^\ast &\bO \\ \endpmatrix , \tag 4.23
 $$
 and, owing to (4.13) and (3.17), (3.18),
 $$
 \gQ^+\cdot1=\mu\bI+\eps\,{(1-\mu)\eta\over\zeta}\,q^{-2r+1+\eps}
 \pmatrix \bI&\bO \\ Z^\ast &\bO \\ \endpmatrix +(1-\mu)
 \pmatrix \bO&\bO \\ V^\ast &\bI \\ \endpmatrix , \tag 4.24
 $$
 (here the dimensions of the zero and unit blocks vary and are determined
 implicitly). Thus, in virtue of (4.2), one can evaluate $M\cdot1$.
 Notice that $M_{jj}\cdot1=\xi_+$, for $j\ge r$. One expects that, moreover,
 $M_{jj}\cdot1=1$, for $-r<j<r$, and $M_{jj}\cdot1=\xi_+^{-1}$, for
 $j\le-r$. It is straightforward to check that this is actually the case
 provided
 $$
 \mu=(1-\xi_0)/(\xi_+-\xi_0) , \quad
 \xi_0=q^{N-2r+1} . \tag 4.25
 $$
 Here we have accounted the equations (4.16b,c) to express $\eta$ and
 $\zeta$.

 \bigpagebreak
 \flushpar {\bf 5. $\CCah$ as a left $\tAd$-module}
 \medpagebreak

 As a first step we shall rewrite the relation (4.3).

 \proclaim{Lemma 5.1}
 The relation adjoint to (4.3) is equivalent to
$$
[\bI-(\bI-\gP^-_2)R_{12}\gZ^\ast_2R_{12}^{-1}]\gM_1R_{21}^{-1}
(\gP^-_2+\gZ^\ast_2)=\gP^-_2\gM_1R_{21}^{-1}
\tag 5.1 $$
Furthermore, the matrix $\bI-(\bI-\gP^-_2)R_{12}\gZ^\ast_2R_{12}^{-1}$
is invertible and it holds true that
 $$\align
[\bI-(\bI-\gP^-_2)R_{12}\gZ^\ast_2R_{12}^{-1}]^{-1}\gP^-_2 = &
R_{12}(\gP^-_2+\gZ^\ast_2)R_{21} \\
& -(q-q^{-1})(\gP^-_1+\gZ^\ast_1)R_{12}P_{12}(\gP^-_1+\gZ^\ast_1) \\
& +(q-q^{-1})\eps q^{-2r+1+\eps}R_{12}(\gP^-_2+\gZ^\ast_2)K_{21}P_{12}
(\gP^-_2+\gZ^\ast_2)R_{21}.\\
 \tag 5.2
 \endalign$$
 \endproclaim

 \demo{Proof}
 The adjoint to (4.3) means that
 $$
(\bI-\gP^-_2-\gZ^\ast_2)R_{12}^{-1}\gM_1R_{21}^{-1}(\gP^-_2+\gZ^\ast_2)=0
 $$
or, equivalently,
$$
(\bI-\gP^-_2-\gZ^\ast_2)R_{12}^{-1}\gM_1R_{21}^{-1}\gZ^\ast_2=
-(\bI-\gP^-_2-\gZ^\ast_2)R_{12}^{-1}\gM_1R_{21}^{-1}\gP^-_2.
 $$
Multiplying this equality from the left by $(\bI-\gP^-_2)R_{12}$ or,
in the opposite direction, by $(\bI-\gP^-_2)R_{12}^{-1}$ and noting that
(cf. (2.22))
$$
(\bI-\gP^-_2)R_{12}^{\pm 1}(\bI-\gP^-_2)R_{12}^{\mp 1} =\bI-\gP^-_2
$$
and $(\bI-\gP^-_2)R_{21}^{-1}\gZ^\ast_2=R_{21}^{-1}\gZ^\ast_2$ one finds that
it is equivalent to
$$
[\bI-(\bI-\gP^-_2)R_{12}\gZ^\ast_2R_{12}^{-1}]\gM_1R_{21}^{-1}
\gZ^\ast_2=-(\bI-\gP^-_2)(\bI-R_{12}\gZ^\ast_2R_{12}^{-1})
\gM_1R_{21}^{-1}\gP^-_2.
$$
{}From the last equality one can easily derive (5.1).

It is easy to show that
$[(\bI-\gP^-_2)R_{12}\gZ^\ast_2R_{12}^{-1}]^4=0$ and so
$\bI-(\bI-\gP^-_2)R_{12}\gZ^\ast_2R_{12}^{-1}$ is invertible. This is
clear also from the following calculation. Using $(\gZ^\ast)^2=0$ one
finds that
$$
[(\bI-\gP^-_2)R_{12}\gZ^\ast_2R_{12}^{-1}]^k\gP^-_2=(-1)^{k-1}
(\bI-\gP^-_2)(R_{12}\gZ^\ast_2R_{12}^{-1}\gP^-_2)^k
\tag 5.3 $$
and consequently
$$\align
[\bI-(\bI-\gP^-_2)R_{12}\gZ^\ast_2R_{12}^{-1}]^{-1}\gP^-_2 & =
\gP^-_2 +(\bI-\gP^-_2)R_{12}\gZ^\ast_2R_{12}^{-1}\gP^-_2
(\bI+R_{12}\gZ^\ast_2R_{12}^{-1}\gP^-_2)^{-1} \\
& = R_{12}(\gP^-_2+\gZ^\ast_2)R_{12}^{-1}\gP^-_2
(\bI+R_{12}\gZ^\ast_2R_{12}^{-1}\gP^-_2)^{-1}. \\
\tag 5.4
\endalign $$
Applying (5.3) in the reversed direction to the expression
$(\bI+R_{12}\gZ^\ast_2R_{12}^{-1}\gP^-_2)^{-1}$ occuring after the last
equality sign in (5.4) one derives that
$$
[\bI-(\bI-\gP^-_2)R_{12}\gZ^\ast_2R_{12}^{-1}]^{-1}\gP^-_2
=R_{12}(\gP^-_2+\gZ^\ast_2)R_{12}^{-1}
[\bI-(\bI-\gP^-_2)R_{12}\gZ^\ast_2R_{12}^{-1}]^{-1}\gP^-_2.
\tag 5.5 $$
Replace $R_{12}^{-1}$ standing between the brackets on the RHS of (5.5)
by $R_{21}-\gamma P_{12}+\gamma K_{21}$. The result is (using
$\gP^-_1\gP^-_2R_{12}\gZ^\ast_2=0$ in the second summand)
$$ \align
& R_{12}(\gP^-_2+\gZ^\ast_2)R_{21}-
\gamma R_{12}(\gP^-_2+\gZ^\ast_2)P_{12}
(\bI+R_{12}\gZ^\ast_2R_{12}^{-1})\gP^-_2 \\
& +\gamma R_{12}(\gP^-_2+\gZ^\ast_2)K_{21}
[\bI-(\bI-\gP^-_2)R_{12}\gZ^\ast_2R_{12}^{-1}]^{-1}\gP^-_2.\\
\endalign $$
Owing to (3.15) the second summands equals
$-\gamma(\gP^-_1+\gZ^\ast_1)R_{12}P_{12}(\gP^-_1+\gZ^\ast_1)$.
To complete the verification of (5.2) it remains to show that
$$
K_{21}[\bI-(\bI-\gP^-_2)R_{12}\gZ^\ast_2R_{12}^{-1}]^{-1}\gP^-_2=
pK_{21}P_{12}(\gP^-_2+\gZ^\ast_2)R_{21}.
$$
In virtue of (5.4) this is the same as
$$
K_{21}R_{12}(\gP^-_2+\gZ^\ast_2)R_{12}^{-1}\gP^-_2=
pK_{21}P_{12}(\gP^-_2+\gZ^\ast_2)R_{21}
(\bI+R_{12}\gZ^\ast_2R_{12}^{-1}\gP^-_2).
$$
Since $pK_{21}P_{12}\gP^-_2R_{21}=K_{21}\gP^-_2$ (cf. (2.24)) the last
equality can be rewritten as
$$
K_{21}(\bI-\gP^-_2)R_{12}(\gP^-_2+\gZ^\ast_2)R_{12}^{-1}\gP^-_2=
pK_{21}P_{12}\gZ^\ast_2R_{21}
(\bI+R_{12}\gZ^\ast_2R_{12}^{-1})\gP^-_2.
\tag 5.6 $$
Using (2.25), $(\bI-\gP^-)\gZ^\ast=\gZ^\ast$ and (2.24) one finds that the
LHS of (5.6) equals
$$ \align
& pK_{21}P_{12}\gZ^\ast_2R_{12}^{-1}\gP^-_2-
\gamma K_{21}P_{12}\gP^-_1\gZ^\ast_2R_{12}^{-1}\gP^-_2 \\
& = pK_{21}P_{12}\gZ^\ast_2R_{12}^{-1}\gP^-_2 -
\gamma pK_{21}\gP^-_1R_{12}\gZ^\ast_2R_{12}^{-1}\gP^-_2. \\
\endalign $$
Substitute $R_{12}^{-1}+\gamma P_{12}-\gamma K_{21}$ for $R_{21}$ on the RHS
of (5.6) to get
$$
 pK_{21}P_{12}\gZ^\ast_2R_{12}^{-1}\gP^-_2
+\gamma pK_{21}\gZ^\ast_1(\bI+R_{12}\gZ^\ast_2R_{12}^{-1})\gP^-_2.
$$
The second summand in the last expression can be simplified using (3.16) as
$$ \align
& \gamma pK_{21}(\gZ^\ast_1-\gP^-_1R_{12}\gZ^\ast_2R_{12}^{-1}-
\gZ^\ast_1R_{12}\gP^-_2R_{12}^{-1})\gP^-_2 \\
& =-\gamma pK_{21}\gP^-_1R_{12}\gZ^\ast_2R_{12}^{-1}\gP^-_2. \\
\endalign $$
This verifies (5.6) and hence (5.2) as well. \qed
 \enddemo

We shall need the following identity.
 \proclaim{Lemma 5.2}
It holds true that
 $$\align
& [\bI-(\bI-\gP^-_2)R_{12}\gZ^\ast_2R_{12}^{-1}]^{-1}\gP^-_2 R_{23}
[\bI-(\bI-\gP^-_3)R_{13}\gZ^\ast_3R_{13}^{-1}]^{-1}\gP^-_3 R_{32} \\
& = R_{23}[\bI-(\bI-\gP^-_3)R_{13}\gZ^\ast_3R_{13}^{-1}]^{-1}\gP^-_3 R_{32}
[\bI-(\bI-\gP^-_2)R_{12}\gZ^\ast_2R_{12}^{-1}]^{-1}\gP^-_2 . \\
 \endalign $$
 \endproclaim

 \demo{Proof}
Because of Lemma 5.1 and using  the notation introduced in Sec.3 we have
to show that
$$ \align
& (R_{12}\gX_2R_{21}-\gamma \gX_1R_{12}P_{12}\gX_1+
\gamma pR_{12}\gX_2K_{21}P_{12}\gX_2R_{21})R_{23} \\
& \times\, (R_{13}\gX_3R_{31}-\gamma \gX_1R_{13}P_{13}\gX_1+
\gamma pR_{13}\gX_3K_{31}P_{13}\gX_3R_{31})R_{32} \\
& = R_{23} (R_{13}\gX_3R_{31}-\gamma \gX_1R_{13}P_{13}\gX_1+
\gamma pR_{13}\gX_3K_{31}P_{13}\gX_3R_{31}) \\
& \times\, R_{32} (R_{12}\gX_2R_{21}-\gamma \gX_1R_{12}P_{12}\gX_1+
\gamma pR_{12}\gX_2K_{21}P_{12}\gX_2R_{21}) . \\
\tag 5.7
\endalign $$
Unfortunately, I was not able to find other proof than that one based
on a straightforward and rather tedious calculation. Below I confine
myself to listing the intermediate equalities which altogether imply
the identity (5.7) and to giving hints just for several less obvious
steps. It holds true that
$$ \align
& R_{12}\gX_2R_{21}R_{23}R_{13}\gX_3R_{31}R_{32}=
R_{23}R_{13}\gX_3R_{31}R_{32}R_{12}\gX_2R_{21},
\tag 5.8 \\
& R_{12}\gX_2K_{21}P_{12}\gX_2R_{21}R_{23}\gX_1R_{13}P_{13}\gX_1R_{32}=0,
\tag 5.9a \\
& R_{23}\gX_1R_{13}P_{13}\gX_1R_{32}R_{12}\gX_2K_{21}P_{12}\gX_2R_{21}=0,
\tag 5.9b \\
& \gX_1R_{12}P_{12}\gX_1R_{23}R_{13}\gX_3K_{31}P_{13}\gX_3R_{31}R_{32}=0,
\tag 5.10a \\
& R_{23}R_{13}\gX_3K_{31}P_{13}\gX_3R_{31}R_{32}\gX_1R_{12}P_{12}\gX_1=0,
\tag 5.10b \\
& \gX_1R_{12}P_{12}\gX_1R_{23}R_{13}\gX_3R_{31}R_{32}=
R_{23}R_{13}\gX_3R_{31}R_{32}\gX_1R_{12}P_{12}\gX_1,
\tag 5.11 \\
& R_{12}\gX_2R_{21}R_{23}\gX_1R_{13}P_{13}\gX_1R_{32}-
R_{23}\gX_1R_{13}P_{13}\gX_1R_{32}R_{12}\gX_2R_{21} \\
& \ =\gamma \gX_1R_{12}R_{13}\gX_2R_{23}\gX_3
(R_{21}P_{12}P_{13}-R_{32}P_{13}P_{12}),
\tag 5.12a \\
& \gX_1R_{12}P_{12}\gX_1R_{23}\gX_1R_{13}P_{13}\gX_1R_{32}=
\gX_1R_{12}R_{13}\gX_2R_{23}\gX_3R_{21}P_{12}P_{13},
\tag 5.12b \\
& R_{23}\gX_1R_{13}P_{13}\gX_1R_{32}\gX_1R_{12}P_{12}\gX_1=
\gX_1R_{12}R_{13}\gX_2R_{23}\gX_3R_{32}P_{13}P_{12},
\tag 5.12c \\
& R_{12}\gX_2K_{21}P_{12}\gX_2R_{21}R_{23}R_{13}\gX_3R_{31}R_{32}=
R_{23}R_{13}\gX_3R_{31}R_{32}R_{12}\gX_2K_{21}P_{12}\gX_2R_{21}.
\tag 5.13 \\
\endalign $$

To obtain few last equalities one can use that $\gY^2=\gY$ and so
$$ \align
K_{21}P_{12}\gX_2R_{21}R_{23}R_{13}\gX_3K_{31}P_{13} & =
p^{-2}K_{21}\gY_1R_{23}\gY_1K_{13} \\
& = p^{-1}K_{21}R_{23}R_{13}\gX_3K_{31}P_{13} \\
& = p^{-1}K_{21}\gX_3K_{31}P_{13} . \\
\endalign $$
Furthermore, observe that
$$
\gP^-_2\gP^-_3K_{21}K_{31}=\gP^-_2\gP^-_3P_{23}K_{31} .
\tag 5.14 $$
In fact, one can readily check that even the equality
$K_{21}K_{31}=P_{23}K_{31}$
holds true but (5.14) can be verified directly as follows
$$ \align
\gP^-_2\gP^-_3K_{21}K_{31} & =
\gamma^{-1}\gP^-_2\gP^-_3(R_{12}^{-1}-R_{21})K_{31} \\
& = \gamma^{-1}\gP^-_2\gP^-_3(R_{32}-R_{23}^{-1})K_{31} \\
& = \gP^-_2\gP^-_3P_{23}K_{31} .
\endalign $$

So finally we have
$$ \align
& R_{12}\gX_2R_{21}R_{23}R_{13}\gX_3K_{31}P_{13}\gX_3R_{31}R_{32}-
R_{23}R_{13}\gX_3K_{31}P_{13}\gX_3R_{31}R_{32}R_{12}\gX_2R_{21} \\
& \ =-\gamma P_{23}R_{13}\gX_3K_{31}\gX_2\gX_1R_{13}R_{12}P_{13}+
\gamma R_{23}R_{13}\gX_3K_{31}\gX_2\gX_1R_{13}P_{12}P_{13},
\tag 5.15a \\
& pR_{12}\gX_2K_{21}P_{12}\gX_2R_{21}R_{23}R_{13}
\gX_3K_{31}P_{13}\gX_3R_{31}R_{32} \\
& \ =P_{23}R_{13}\gX_3K_{31}\gX_2\gX_1R_{13}R_{12}P_{13},
\tag 5.15b \\
& pR_{23}R_{13}\gX_3K_{31}P_{13}\gX_3R_{31}R_{32}R_{12}\gX_2K_{21}P_{12}
\gX_2R_{21} \\
& \ =R_{23}R_{13}\gX_3K_{31}\gX_2\gX_1R_{13}P_{12}P_{13} .
\tag 5.15c \\
\endalign $$
The equalities (5.8), (5.9a,b), (5.10a,b), (5.11), (5.12a,b,c), (5.13) and
(5.15a,b,c) imply (5.7). \qed
\enddemo

Now we can equip $\CCah$ with the structure of a left $\tAd$-module.
First we should specify the value of $M\cdot 1$. The result is suggested
by (4.2), (4.23) and (4.24). Accepting (4.25), the value will depend
just on one parameter $\xi\equiv\xi_+\not=0$. To extend the action to
the whole algebra $\CCah$ one can use Lemma 5.1. Thus we postulate
$$ \align
& M\cdot 1= \bI+(\xi^{-1}-1)(\gP^-+\gZ^\ast)+(\xi-1)(\gP^++\gV^\ast),
\tag 5.16 \\
& M_1R_{21}^{-1}\cdot(\gP^-_2+\gZ^\ast_2)f=
[\bI-(\bI-\gP^-_2)R_{12}\gZ^\ast_2R_{12}^{-1}]^{-1}\gP^-_2M_1R_{21}^{-1}
\cdot f,\quad\forall f\in\CCah .
\tag 5.17 \\
\endalign $$
Before verifying that (5.16), (5.17) really define a left module
structure let us show
\proclaim{Lemma 5.3}
The relation (5.17) implies
$$
M_1R_{21}^{-1}\cdot(\gP^+_2+\gV^\ast_2)f=
\gP^+_2[\bI-R_{12}\gV^\ast_2R_{12}^{-1}(\bI-\gP^+_2)]^{-1}M_1R_{21}^{-1}
\cdot f,\quad\forall f\in\CCah .
\tag 5.18 $$
The matrix $\bI-R_{12}\gV^\ast_2R_{12}^{-1}(\bI-\gP^+_2)$ is invertible
and it holds true that
$$ \align
\gP^+_2[\bI-R_{12}\gV^\ast_2R_{12}^{-1}(\bI-\gP^+_2)]^{-1}= &
R_{21}^{-1}(\gP^+_2+\gV^\ast_2)R_{12}^{-1} \\
& +(q-q^{-1})(\gP^+_1+\gV^\ast_1)R_{21}^{-1}P_{12}(\gP^+_1+\gV^\ast_1) \\
& -(q-q^{-1})\eps q^{2r-1-\eps}R_{21}^{-1}(\gP^+_2+\gV^\ast_2)K_{21}P_{12}
(\gP^+_2+\gV^\ast_2)R_{12}^{-1} .\\
 \tag 5.19 \endalign
$$
\endproclaim

\demo{Proof}
The matrix $\bI-R_{12}\gV^\ast_2R_{12}^{-1}(\bI-\gP^+_2)$ is again
invertible since, as one can easily check, the matrix
$R_{12}\gV^\ast_2R_{12}^{-1}(\bI-\gP^+_2)$ is nilpotent. The equality
(5.19) can be verified in a way quite similar as it has been done
for the relation (5.2) given in Lemma 5.1. We omit the details.

Let us turn to the equality (5.18). One can start from (cf. (5.17))
$$
pM_3R_{23}^{-1}\cdot\gX_2fR_{13}^{-1}K_{21}
=p [\bI-(\bI-\gP^-_2)R_{32}\gZ^\ast_2R_{32}^{-1}]^{-1}\gP^-_2
M_3\cdot fR_{23}^{-1}R_{13}^{-1}K_{21} .
\tag 5.20 $$
The LHS equals
$$
M_3\cdot R_{23}^{-1}R_{13}^{-1}R_{12}^{-1}\gY_1P_{12}K_{21}f=
R_{12}^{-1}M_3\cdot R_{13}^{-1}\gY_1R_{13}P_{12}K_{21}
$$
and so it is possible to rewrite (5.20) as
$$
M_3\cdot R_{23}^{-1}\gY_2R_{23}K_{21}f=
pP_{12}R_{12}[\bI-(\bI-\gP^-_2)R_{32}\gZ^\ast_2R_{32}^{-1}]^{-1}\gP^-_2
K_{21}M_3\cdot f .
\tag 5.21 $$
Next we take into account the identity (5.2). The following
equalities are straightforward to verify:
$$ \align
& pP_{12}R_{12}R_{32}\gX_2R_{23}K_{21}=R_{23}^{-1}\gY_2R_{32}^{-1}K_{21},\\
& pP_{12}R_{12}\gX_3R_{32}P_{23}\gX_3K_{21}=
p\gX_3P_{12}K_{31}P_{23}\gX_3, \\
& p^2P_{12}R_{12}R_{32}\gX_2K_{23}P_{23}\gX_2R_{23}K_{21}=
\gY_3R_{23}^{-1}\gY_2K_{31}K_{21} . \\
\endalign $$
Now it is enough to take the trace in the first factor of (5.21) and to use
$$
\tr_1\,K_{12}=\bI_2,\ \tr_1(P_{12}K_{31})=K_{32},\
\tr_1(K_{31}K_{21})=P_{23} .
\tag 5.22 $$
Observe also that
$$
p\gX_3K_{32}P_{23}\gX_3=p^{-1}R_{23}^{-1}\gY_2K_{23}P_{23}\gY_2R_{32}^{-1} .
\tag 5.23 $$
This way one obtains
$$
M_1R_{21}^{-1}\cdot\gY_2f=
(R_{21}^{-1}\gY_2R_{12}^{-1}+\gamma\gY_1R_{21}^{-1}P_{12}\gY_1-
\gamma p^{-1}R_{21}^{-1}\gY_2K_{21}P_{12}\gY_2R_{12}^{-1})M_1R_{21}^{-1}
\cdot f .
$$
The relation (5.18) then follows from (5.19). \qed
\enddemo

Now we can state
\proclaim{Proposition 5.4}
The relations (5.16), (5.17) define unambiguously on $\CCah$ the structure
of a left $\tAd$-module depending on one parameter $\xi\not=0$.
\endproclaim

\demo{Proof}
First we have to show that (5.16), (5.17) define correctly a linear mapping
$\CCah\to\Mat(N)\otimes\CCah: f\mapsto M\cdot f$. Let $\bCCah$ be the free
algebra generated by the entries of the matrix
$\bar Z^\ast=(\bar Z^\ast_{jk})=(\bza_{kj})$, with $-r<j,\ k\le -r$.
Hence $\CCah$ is obtained
from $\bCCah$ by means of factorization by the two-sided ideal generated
by the
relation adjoint to (3.15) and $z^\ast_{kj}$ are the factor images of
$\bza_{kj}$. In the obvious notation, the matrix $\bgZa$ is obtained
from $\gZ^\ast$ when replacing all elements $z^\ast_{kj}$ by $\bza_{kj}$.
Set (in this proof) $\bgX:=\gP^-+\bgZa$ and
$$
\bX_{12}:=[\bI-(\bI-\gP^-_2)R_{12}\bgZa_2R_{12}^{-1}]^{-1}\gP^-_2 .
$$
It is clear that when replacing everywhere in (5.16), (5.17) the elements
$z^\ast_{kj}$ by $\bza_{kj}$ one obtains a well defined linear mapping
$$
\bCCah\to\Mat(N)\otimes\bCCah: f\mapsto M\cdot f .
\tag 5.24 $$
A straightforward calculation gives
$$ \align
& M_1\cdot R_{21}^{-1}R_{31}^{-1}(\bgX_2R_{23}\bgX_3R_{32}-
R_{23}\bgX_3R_{32}\bgX_2)f \\
& = (\bX_{12}R_{23}\bX_{13}R_{32}-R_{23}\bX_{13}R_{32}\bX_{12})M_1\cdot
R_{21}^{-1}R_{31}^{-1}f,\quad\forall f\in\bCCah .
\endalign $$
This means that the linear mapping (5.24) factorizes from $\bCCah$ to
$\CCah$ if and only if the factor-image of the matrix
$\bX_{12}R_{23}\bX_{13}R_{32}-R_{23}\bX_{13}R_{32}\bX_{12}$ vanishes.
But this is the content of Lemma 5.2.

To show that $\CCah$ is really a left $\tAd$-module we have to verify
the equality
$$
M_2R_{12}^{-1}M_1R_{21}^{-1}\cdot1=R_{12}^{-1}M_1R_{21}^{-1}M_2\cdot1
\tag 5.25 $$
and the implication (some parentheses appear here and in what follows
only for graphical reasons)
$$ \align
& (M_2R_{12}^{-1}M_1R_{21}^{-1})R_{31}^{-1}R_{32}^{-1}\cdot f=
(R_{12}^{-1}M_1R_{21}^{-1}M_2)R_{31}^{-1}R_{32}^{-1}\cdot f\Longrightarrow \\
& (M_2R_{12}^{-1}M_1R_{21}^{-1})R_{31}^{-1}R_{32}^{-1}\cdot \gZ^\ast_3f=
(R_{12}^{-1}M_1R_{21}^{-1}M_2)R_{31}^{-1}R_{32}^{-1}\cdot \gZ^\ast_3f ,\quad
\forall f\in\CCah ,
\tag 5.26 \endalign
$$
since then (5.25) and (5.26) jointly imply
$$
M_2R_{12}^{-1}M_1R_{21}^{-1}\cdot f=R_{12}^{-1}M_1R_{21}^{-1}M_2\cdot f,
\quad\forall f\in\CCah .
$$
In this proof we set
$$\align
& X_{12}:=[\bI-(\bI-\gP^-_2)R_{12}\gZ^\ast_2R_{12}^{-1}]^{-1}\gP^-_2 ,\\
& Y_{12}:=\gP^+_2[\bI-R_{12}\gV^\ast_2R_{12}^{-1}(\bI-\gP^+_2)]^{-1} .\\
\endalign $$

{\it Verification of (5.25)}.
Using (5.16), (5.17) and (5.18) one finds immediately that
$$ \align
& M_2R_{12}^{-1}M_1R_{21}^{-1}\cdot1=
\big(\bI+(\xi^{-1}-1)X_{21}+(\xi-1)Y_{21}\big)
\big(\bI+(\xi^{-1}-1)\gX_2+(\xi-1)\gY_2\big)R_{12}^{-1}R_{21}^{-1} ,\\
& R_{12}^{-1}M_1R_{21}^{-1}M_2\cdot1=R_{12}^{-1}
\big(\bI+(\xi^{-1}-1)X_{12}+(\xi-1)Y_{12}\big)
\big(\bI+(\xi^{-1}-1)\gX_1+(\xi-1)\gY_1\big)R_{21}^{-1} . \\
\endalign $$
This means that (5.25) holds for every $\xi\not=0$ if and only if the
following five equation are satisfied:
$$ \align
& X_{21}\gX_2R_{12}^{-1}=R_{12}^{-1}X_{12}\gX_1 ,\tag 5.27 \\
& X_{21}(\bI-\gX_2-\gY_2)R_{12}^{-1}+(\bI-X_{21}-Y_{21})\gX_2R_{12}^{-1} \\
& \ =R_{12}^{-1}X_{12}(\bI-\gX_1-\gY_1)+R_{12}^{-1}
(\bI-X_{12}-Y_{12})\gX_1 ,\tag 5.28 \\
& Y_{21}\gX_2R_{12}^{-1}+(\bI-X_{21}-Y_{21})(\bI-\gX_2-\gY_2)R_{12}^{-1}+
X_{21}\gY_2R_{12}^{-1} \\
& \ =R_{12}^{-1}Y_{12}\gX_1+R_{12}^{-1}(\bI-X_{12}-Y_{12})
(\bI-\gX_1-\gY_1)+R_{12}^{-1}X_{12}\gY_1 ,\tag 5.29 \\
& Y_{21}(\bI-\gX_2-\gY_2)R_{12}^{-1}+(\bI-X_{21}-Y_{21})\gY_2R_{12}^{-1} \\
& \ =R_{12}^{-1}Y_{12}(\bI-\gX_1-\gY_1)+
R_{12}^{-1}(\bI-X_{12}-Y_{12})\gY_1 , \tag 5.30 \\
& Y_{21}\gY_2R_{12}^{-1}=R_{12}^{-1}Y_{12}\gY_1 . \tag 5.31 \\
\endalign $$
The equations (5.27-31) are not independent. In fact, it is enough to
verify (5.27), (5.31) and the following two equations:
$$ \align
& (X_{21}+\gX_2-X_{21}\gY_2-Y_{21}\gX_2)R_{12}^{-1}=
R_{12}^{-1}(X_{12}+\gX_1-X_{12}\gY_1-Y_{12}\gX_1) ,\tag 5.32 \\
& (Y_{21}+\gY_2-Y_{21}\gX_2-X_{21}\gY_2)R_{12}^{-1}=
R_{12}^{-1}(Y_{12}+\gY_1-Y_{12}\gX_1-X_{12}\gY_1) .\tag 5.33 \\
\endalign $$
This can be done by a straightforward calculation using (5.2) and (5.19).
We omit the details.

{\it Verification of (5.26)}.
Let us first show that
$$
(\bI-\gX_3)R_{23}^{-1}R_{13}^{-1}(M_2R_{12}^{-1}M_1R_{21}^{-1})
R_{31}^{-1}R_{32}^{-1}\cdot\gX_3f=0 .
\tag 5.34 $$
Actually, (5.17) means that (check the proof of Lemma 5.1)
$$
M_1R_{21}^{-1}\cdot\gX_2f=R_{12}\gX_2R_{12}^{-1}M_1R_{21}^{-1}\cdot\gX_2f
$$
and so
$$ \align
(M_2R_{12}^{-1}M_1R_{21}^{-1})R_{31}^{-1}R_{32}^{-1}\cdot\gX_3f & =
M_2R_{12}^{-1}R_{32}^{-1}M_1R_{31}^{-1}\cdot\gX_3R_{21}^{-1}f \\
& =M_2\cdot(R_{13}R_{32}^{-1}R_{12}^{-1}\gX_3R_{13}^{-1}M_1R_{31}^{-1}\cdot
\gX_3R_{21}^{-1}f) \\
& =R_{13}R_{23}\gX_3R_{23}^{-1}M_2\cdot
(R_{32}^{-1}\gX_3R_{12}^{-1}R_{13}^{-1}M_1R_{31}^{-1}
\cdot\gX_3R_{21}^{-1}f) . \\
\endalign $$
Now multiply the both sides by $(\bI-\gX_3)R_{23}^{-1}R_{13}^{-1}$
from the left and recall that $(\bI-\gX)\gX=0$.

Further one can proceed similarly as in the proof of Lemma 5.1. The
relation (5.34) is equivalent to
$$ \align
& (\bI-\gP^-_3-\gZ^\ast_3)R_{23}^{-1}R_{13}^{-1}(M_2R_{12}^{-1}M_1
R_{21}^{-1})R_{31}^{-1}R_{32}^{-1}\cdot\gZ^\ast_3f \\
& =-(\bI-\gP^-_3-\gZ^\ast_3)R_{23}^{-1}R_{13}^{-1}(M_2R_{12}^{-1}M_1
R_{21}^{-1})R_{31}^{-1}R_{32}^{-1}\gP^-_3\cdot f . \\
\endalign $$
Multiplying this equation by $(\bI-\gP^-_3)R_{13}R_{23}$ from the left
and taking into account that
$$ \align
& (\bI-\gP^-_3)R_{13}R_{23}(\bI-\gP^-_3)R_{23}^{-1}R_{13}^{-1}=
\bI-\gP^-_3 , \\
& (\bI-\gP^-_3)R_{31}^{-1}R_{32}^{-1}\gZ^\ast_3=
R_{31}^{-1}R_{32}^{-1}\gZ^\ast_3 ,\\
\endalign $$
one arrives at
$$ \align
& [\bI-(\bI-\gP^-_3)R_{13}R_{23}\gZ^\ast_3R_{23}^{-1}R_{13}^{-1}]
(M_2R_{12}^{-1}M_1R_{21}^{-1})R_{31}^{-1}R_{32}^{-1}\cdot\gZ^\ast_3f \\
& =-(\bI-\gP^-_3)(\bI-R_{13}R_{23}\gZ^\ast_3R_{23}^{-1}R_{13}^{-1})
(M_2R_{12}^{-1}M_1R_{21}^{-1})R_{31}^{-1}R_{32}^{-1}\gP^-_3\cdot f .
\tag 5.35 \\
\endalign $$
Quite similarly one obtains

$$ \align
& [\bI-(\bI-\gP^-_3)R_{13}R_{23}\gZ^\ast_3R_{23}^{-1}R_{13}^{-1}]
(R_{12}^{-1}M_1R_{21}^{-1}M_2)R_{31}^{-1}R_{32}^{-1}\cdot\gZ^\ast_3f \\
& =-(\bI-\gP^-_3)(\bI-R_{13}R_{23}\gZ^\ast_3R_{23}^{-1}R_{13}^{-1})
(R_{12}^{-1}M_1R_{21}^{-1}M_2)R_{31}^{-1}R_{32}^{-1}\gP^-_3\cdot f .
\tag 5.36 \\
\endalign $$
The right hand sides of (5.35) and (5.36) are equal by assumption and so,
to complete the proof, it suffices to show that that the matrix
$\bI-(\bI-\gP^-_3)R_{13}R_{23}\gZ^\ast_3R_{23}^{-1}R_{13}^{-1}$
is invertible. To this end, observe from (2.3) that the R matrix $R_{12}$
is lower triangular provided the lexicographical ordering of the basis
 vectors in the tensor product has been chosen
($R_{jk,st}=0$ for $j<s$, and $R_{jk,jt}=\delta_{kt}R_{jk,jt}$).
Hence in the tensor product $\BC^N\otimes\BC^N\otimes\BC^N$, the matrix
$\bI-\gP^-_3$ is diagonal, the matrices $R_{13}^{\pm1},\ R_{23}^{\pm1}$
are lower triangular and the matrix $\gZ^\ast_3$ is lower triangular
with vanishing diagonal. Consequently, the matrix
$(\bI-\gP^-_3)R_{13}R_{23}\gZ^\ast_3R_{23}^{-1}R_{13}^{-1}$
is lower triangular with vanishing diagonal and hence nilpotent. \qed
\enddemo

\newpage
 \bigpagebreak
 \flushpar {\bf 6. Irreducible representations}
 \medpagebreak

Set $\xi=q^{-2\sigma}$. Observe that $M\cdot1$ is a lower triangular matrix
(cf. (5.16)). Since $M=\Lambda^\ast\Lambda$ and in view of the form of the
isomorphism (2.29), (2.30), this means that $X^-_j\cdot1=0,\ \forall j$, and
so $1$ is a lowest weight vector with a lowest weight $\lambda$. Denote by
$\CMs$ the cyclic $\tAd$-submodule in $\CCah$ with the cyclic vector 1.
Owing to the isomorphism (2.29), (2.30), $\CMs$ can be also regarded as a
$\Uh$-module. It remains to determine the conditions
implying that the module $\CMs$ is finite-dimensional. This means to
determine the parameters $\sigma$. Recall that by the results due to Rosso
\cite{\Rosso}, $\CMs$ is determined unambiguously, up to equivalence, by the
lowest weight  $\lambda$.

In more detail,
 $$\align
 M_{jj}\cdot1 & = q^{2\sigma}\ ,\text{ for }\ j\le-r, \\
 & =1\ \quad,\qquad -r<j<r, \tag 6.1\\
 & =q^{-2\sigma},\ \qquad j\ge r. \\
 \endalign $$
Consequently
$$\align
 q^{H_j}\cdot1 & =q^{-\sigma}\ ,\text{ for }\ j=r,\tag 6.2\\
 & =1\ \quad,\qquad j\not=r,\\
 \endalign$$
 in the case of $B_l$;
 $$\left.\aligned
 q^{H_j}\cdot1 & =q^{-2\sigma}\ ,\text{ for }\ j=1,\\
 & =1\quad\quad,\qquad j>1,\\
 \endaligned\right\}\text{ if }r=\half\, ;\tag 6.3a
 $$
 $$\left.\aligned
 q^{H_j}\cdot1 & =q^{-\sigma}\ ,\text{ for }\ j=r+\half,\\
 & =1\quad\quad,\qquad j\not=r+\half,\\
 \endaligned\right\}\text{ if }r>\half\, ;\tag 6.3b
 $$
 in the case of $C_l$;
 $$\left.\aligned
 q^{H_j}\cdot1 & =q^{-2\sigma}\ ,\text{ for }\ j=1,\\
 & =1\quad\quad,\qquad j>1,\\
 \endaligned\right\}\text{ if }r=\half\, ;\tag 6.4a
 $$
 $$\left.\aligned
 q^{H_j}\cdot1 & =q^{-\sigma}\ ,\text{ for }\ j=1,\\
 & =q^{-\sigma}\ ,\ \qquad j=2,\\
 & =1\ \quad,\qquad j>2,\\
 \endaligned\right\}\text{ if }r=\thalf\, ;\tag 6.4b
 $$
 $$\left.\aligned
 q^{H_j}\cdot1 & =q^{-\sigma}\ ,\text{ for }\ j=r+\half,\\
 & =1\quad\quad,\qquad j\not=r+\half,\\
 \endaligned\right\}\text{ if }r>\thalf\, ;\tag 6.4c
 $$
 in the case of $D_l$.

Let $\{\omega_1,\dots,\omega_l\}\subset\gh^\ast$ be the set of fundamental
weights corresponding to the set of simple roots
$\Pi=\{\alpha_1,\dots,\alpha_l\}$. Its values on the basis
$\{H_1,\dots,H_l\}$ of the Cartan algebra $\gh$ are given by
$$
\omega_j(H_k)=\half\,\langle\alpha_k,\alpha_k\rangle\,\delta_{jk}.
$$
By our choice, specified in Sec. 2, we have
$$\aligned
& \text{the case } B_l:\ \langle\alpha_1,\alpha_1\rangle=1,\
\langle\alpha_j,\alpha_j\rangle=2,\ j>1;\\
& \text{the case } C_l:\ \langle\alpha_1,\alpha_1\rangle=4,\
\langle\alpha_j,\alpha_j\rangle=2,\ j>1;\\
& \text{the case } D_l:\
\langle\alpha_j,\alpha_j\rangle=2,\ \forall j.\\
\endaligned $$

Thus we have arrived at
\proclaim{Proposition 6.1}
$\CMs$ is a finite-dimensional irreducible $\Uh$-module provided
 $$
 \sigma\in\left\{ \aligned
 & \half\BZ_+\ ,\text{ for }\ r=1;\text{ then }\ \lambda=-2\sigma\omega_1,\\
 & \ \ \BZ_+\ ,\qquad r>1;\qquad\quad \lambda=-\sigma\omega_r,
 \endaligned \right. \tag 6.5
 $$
 in the case of $B_l$;
 $$\aligned
 \sigma & \in \BZ_+\ ,\text{ for }\ \forall r;
\text{ then }\ \lambda=-\sigma\omega_{r+\half},\\
 \endaligned\tag 6.6
 $$
 in the case of $C_l$;
 $$
\sigma\in\left\{\aligned
 & \half\BZ_+\ ,\text{ for }\ r=\half;\text{ then }\
    \lambda=-2\sigma\omega_1,\\
 & \ \ \BZ_+\ ,\qquad r=\thalf;\qquad\quad
     \lambda=-\sigma(\omega_1+\omega_2),\\
 & \ \ \BZ_+\ ,\qquad r>\thalf;\qquad\quad
     \lambda=-\sigma\omega_{r+\half},
\endaligned \right. \tag 6.7
 $$
in the case of $D_l$.
\endproclaim

 \bigpagebreak
 \flushpar {\bf 7. Concluding remarks}
 \medpagebreak

 Though the main goal -- the explicit derivation of antiholomorphic
representations in the Borel-Weil spirit has been reached
 the presented results are not fully satisfactory. This concerns Section 4.
The weakest point is that
 the commutation relations (4.3), (4.4) are apparently over-determined.
 In the ideal case one should extract from them a minimal set of
 commutation relations between $Z_1$ and $Z_2^\ast$ and show the rest to be a
 consequence; particularly this concerns the equality (4.4) or
 equivalently (4.13). However the semiclassical limit yields a formula
 which is not very encouraging in this respect.

 The semiclassical limit $h\to 0$ is based  on the rules
 ($q=\text{e}^{-h}$; {\bf r} is the classical r-matrix)
 $$\align
 & R=\bI-\text{i}h\br+O(h^2)\, ,\\
 & fg-gf=\text{i}h\{f,g\}+O(h^2)\,.
 \endalign$$
 In the classical case one replaces $\xi_0$ by 1, $\xi_\pm$ by $\xi^{\pm1}$,
 $\gQ_{\text{class}}^-$ is formally parameterized in the same way as in the
 quantum case (4.1) and
 $\gQ_{\text{class}}^+=C^0(\gQ_{\text{class}}^-)^t C^0$. So
 $$
 \gM_{\text{class}}^{\ \pm1}=\bI+
 (\xi^{\pm1}-1)\gQ_{\text{class}}^+ +
 (\xi^{\mp1}-1)\gQ_{\text{class}}^-\, .
 $$
 Applying the above rules to (4.3) one finds easily that
 $$
 \{\gM_1\overset\text{\bf .}\to, Z_2\}=
 \pmatrix \bI&Z\endpmatrix_2(\br_{21}\gM_1+\gM_1\br_{12})
 \pmatrix -Z\\ \bI\endpmatrix_2 \,.
 $$
 It is not difficult to obtain also the bracket
 $\{\gM_1^{\ -1}\overset\text{\bf .}\to, Z_2\}$. Now, since
 $\gQ_{\text{class}}^-$ is expressible as a linear combination of the
 matrices $\bI,\ \gM_{\text{class}}$ and $\gM_{\text{class}}^{\ -1}$,
 one can calculate the bracket
 $\{\gQ_1^-\overset\text{\bf .}\to, Z_2\}$ and hence
 $\{Z_1^\ast\overset\text{\bf .}\to, Z_2\}$. There is no need to give the
 final expression but it is rather awkward.

 However explicit calculations in the lowest rank cases, though not presented
 here, support the conjecture that the definition of the algebra $\CC$
 contains no contradiction. To get a more transparent picture,
 one should perhaps attempt to derive the quantum cell $\CC$ as descendant
 from some simpler and more fundamental structure, like the Euclidean space,
 an idea pursued already in the paper \cite{\FRT}. Another approach notable
 in the current literature suggests that one can try to represent part of the
 algebra $\CC$ with the help of quantum differential operators
 \cite{\Japanese}.
 
 \enddocument